\documentclass[12pt,a4paper,leqno]{article}

\usepackage{fullpage,amsmath,amsthm,amsfonts,amscd,graphicx,
latexsym,amssymb,eucal,curves}
\newcommand{\ZZ}{\mathbb{Z}}

\newcommand{\CC}{\mathbb{C}}

\newcommand{\VV}{\mathbb{V}}

\newcommand{\PicF}{{\rm Pic}}
\newcommand{\PicS}{{\rm Pic}}
\newcommand{\JacS}{{\rm Jac}}
\newcommand{\Corr}{{\rm Corr}}

\newcommand{\AlbS}{{\rm Alb}}

\newcommand{\mfS}{{\mathfrak S}}

\newcommand{\mfC}{{\mathfrak C}}
\newcommand{\mrC}{M}

\newcommand{\Spc}{{\rm Spec \,}}

\newcommand{\sms}{\smallsetminus}

\newcommand{\ol}{\overline}

\renewcommand{\div}{{\rm div}}
\newtheorem{Defi}{Definition}[section]

\newtheorem{Rem}[Defi]{Remark}

\newtheorem{Prop}[Defi]{Proposition}
\newtheorem{Lemma}[Defi]{Lemma}
\newtheorem{Cor}[Defi]{Corollary}
\newtheorem{Thm}[Defi]{Theorem}

\begin{document}{\large}
\title{
%Moduli spaces for 
Isotrivially fibred isles in the moduli space of
surfaces of general type}
\author{Martin M{\"o}ller}
\date{}
\maketitle

%\begin{abstract}
%\noindent
\begin{quote}
{\footnotesize {\bf Abstract}.
We complement Catanese' s results on isotrivially fibred
surfaces by completely describing the components
containing an isotrivial surface with monodromy group
$\ZZ/2 \ZZ$. We also give an example
for deformation equivalent isotrivial surfaces with
different monodromy group. 
%\end{abstract}


}

\end{quote}

\section{Introduction}

In \cite{Ca00}, Catanese introduced the notion of 'surfaces
isogenous to a product', partly because the notion of
'isotrivially fibred surface' is not stable under deformation.
For the former he obtains a complete discription of the moduli space
(of surfaces of general type)
by means of fixed loci in products of Teichm{\"u}ller spaces,
while for isotrivially fibred surfaces his description is 
restricted to the sublocus of surfaces having the same
topological type as a given isotrivially fibred surfaces
minus their singular fibres.
\newline
In this paper we complement these results in two ways.
First we describe entirely the connected components $N_2$ of the moduli
space, which contain these isotrivial isles in the simplest
case of monodromy group $G = \ZZ /2 \ZZ$. We use the fact
(Theorem \ref{discrthm}),  that these surfaces are double coverings of
a product of curves, if we impose some restrictions on
the irregularity. This enables us to enumerate the connected
components of $N_2$, and to calculate its dimension and that
of its isotrivially fibred subloci (Corollaries \ref{corcomp} and
\ref{cordim}).
\newline
Secondly  we give an example, that this method does not extend
to bigger monodromy groups, by showing that isotrivially 
fibred surfaces with different monodromy group can occur
in one connected component of the moduli space (Proposition \ref{propneq}).
For this purpose we use the describtion of abelian coverings
given in \cite{FP97}.
\newline
The results of this paper are part of the author's dissertation (\cite{Mo02}),
which describes moduli spaces for different types of fibred surfaces,
by exploiting high irregularity.
\par

\subsubsection*{Acknowledgements}
The author thanks his thesis advisor Frank Herrlich for his
permanently open door to his office.

\subsection*{Notation}
\begin{itemize}
\item We use the complex numbers $\CC$ as our base field throughout.
%\item $\mfC_g$ denotes the moduli functor for smooth curves of
%genus $g$ and $\mrC_g$ its coarse moduli space.
\item By a surface of general type $X$ we think of its canonical
model, i.e.\ a normal surface with $K_X$ ample and at most
rational double points. 
\item $\mfS(\cdot)$ is the functor, which associates to a scheme $T$
the isomorphism classes of families of surfaces of general
type over T. 'Family' always implies proper and flat. We denote by
$N$ its coarse moduli space.
\item A fibration of $X$ is a morphism $X \to B$ to a smooth curve $B$
of genus $b$ with connected fibres. Let $g$ denote the genus of
the fibres. We call $(g,b)$ the type of the fibration. We will
{\em suppose that $g \geq 2$ (fibre condition (FC)) and $b \geq 2$
(base condition (BC))}, so that by
\cite{Be91} and \cite{Ca00} the property of having a type $(g,b)$-fibration
is deformation invariant. This also implies that minimal
models and relatively minimal models coincide. And they admit
such a fibration if and only if the minimal model does.
\item We denote by $\mfS_{g,b}(\cdot)$ the subfunctor of\
$\mfS(\cdot)$ parametrizing $(g,b)$-fibred surfaces. We only
employ this notation when $(BC)$ and $(FC)$ are fulfilled.
\item We use different notations to distinguish between a 
line bundle $L$ on a surface $X
\in \mfS(\Spc \CC)$ and a
family of line bundles ${\cal L}$ on a family of surfaces $X \to T$. 
Given a locally free sheaf $L$ on $X$, we denote
the associated geometric vector bundle by $\VV(L) \to X$.

\end{itemize}

\par

\section{Monodromy and Albanese-image}

Given a fibration $h: X \to B$, we let $B'$ the locus where
the fibres are smooth. Recall that $h$ is said to be {\em isotrivial}
(or {\em of constant moduli}), if the image of the canonical 
mapping $B' \to \mrC_b$ to the moduli space for smooth curves
is just one point.
\newline
Due to (FC) there is a group $G$, the {\em monodromy group},
such that after a galois base change $C_1 \to B$ with
group $G$, {\'e}tale outside the images of degenerate fibres,
we have 
$$ C_1 \times C_2 \cong C_1 \times_B X,$$
where $C_2$ is a fibre of $h$. $G$ acts on both $C_1$ and $C_2$
and we denote the quotient under the diagonal action by
$Q = (C_1 \times C_2)/G$. 
\newline
For symmetry let $B = B_1$, $b=b_1$, $B_2=C_2/G$ and
$g(B_2) = b_2$. We abbreviate $Y = B_1 \times B_2$ 
and call the projections $p_i: Y \to B_i$.
By \cite{Ser96} we know that $q(X) = b_1 + b_2$ and if $b_2 > 0$ the
natural birational map $\varepsilon: X \to Q$ is a morphism, 
i.e.\ $X$ has two fibrations $h_1: X \to B_1$ and $h_2:
X \to B_2$.
\newline
Details for this can be found in \cite{Ser96} and \cite{Ca00}.
\par

\begin{Lemma} \label{AlbLemma}
Let $X \to T \in \mfS_{g,b}$ be a family of
surfaces. Suppose that (at least) one fibre admits an 
isotrivial fibration with $q \geq b+2$.
If $X \to T$ admits a section, the image of the (relative) 
Albanese-map $\alpha$ is a family of
products of smooth curves $B_i$ ($i=1,2$) 
over $T$ of genera $b_1$ and $b_2$. 
\end{Lemma}
\par
{\bf Proof:} Suppose $T = \Spc \CC$ and $X$ isotrivial. The
Albanese variety has dimension $b_1 + b_2$ and its universal
property applied to $h_1 \times h_2$ induces an isogeny
$\AlbS(X) \to \JacS(B_1) \times \JacS(B_2)$
and hence a finite morphism $\alpha(X) \to Y$. 
The image of $Y$ under the inverse isogeny
generates the Albanese variety. As we excluded the cases
$b_1=1$ and $b_2 =1$, this is only possible when
$\alpha(X) = Y$.
\newline
This argument did not make full use of the isotriviality but only
of the existence of both fibration and numerical conditions, 
i.e.\ deformation invariant properties. Hence it also
applies to the other fibres and in the relative setting
by rigidity of a product of curves.
\hfill $\Box$
\par
\begin{Rem}{ \rm
If $q = b_1 + 1$, the Albanese image need not to be a product
of $B_1$ and an elliptic curve. This is the reason, why we exclude
this case, although similar investigations can be made also for
coverings of elliptic fibrations.
}\end{Rem}
\par 

\section{Isotrivial surfaces with monodromy group $G = \ZZ/2\ZZ$}

{\em We restrict ourselves in this section to surfaces in $\mfS_{g,b}^{2}$,
which is the subfunctor of $\mfS_{g,b}$ with $q \geq b_1+2$
parametrizing surfaces, whose (generic) degree of the Albanese
map equals $2$ and where the corresponding components of the moduli space
contain an isotrivially fibred surface, i.e. such that Lemma
\ref{AlbLemma} applies. We shall exclude the case that
$G$ acts freely on $C_1 \times C_2$ , which is dealth with
in \cite{Ca00} (those surfaces will be called {\em isogenous
to a product}). }
\par
Under these conditions Lemma \ref{AlbLemma} ensures for each
$X \in \mfS_{g,b}^{2}(T)$ the existence of fibrations $h_i: X \to B_i$: 
The section required there can be created after an {\'e}tale base
change and these fibrations (though not the Albanese morphism)
descend to $T$. The next theorem will show, that $\alpha$ is not
only generically of degree $2$ but in fact finite.
\par
Recall that a double covering $X \to Y$ of surfaces over $T$ is 
given by data $({L}, D)$, where $D$ is an effective,
flat divisor on $Y$ and ${L}$ is a line bundle on $Y$ with
${L}^{\otimes 2} = {\cal O}_Y(D)$.
\par
\begin{Thm} \label{discrthm}
The surfaces in $\mfS_{g,b}^2(\CC)$ are double coverings of a product
of curves $Y = B_1 \times B_2$ of genus $b_i$, $b_i \geq 2$, ramified
over a curve with at most simple singularities.
\newline
A general surface in $\mfS_{g,b}^2(\CC)$ is smooth, i.e.\ given
$X_0 \in \mfS_{g,b}^2(\CC)$, there is a family $X \to T \in \mfS_{g,b}^2
(T)$ over a $1$-dimensional pointed base $(T,0)$, 
whose fibre over $0$ is $X_0$ and where
all the other fibres are smooth.
\newline
A surface $X$ in $\mfS_{g,b}^2(\CC)$ is isotrivially fibred, if and
only if the branch divisor of $X \to Y$ is composed of horizontal
and vertical curves. 
\end{Thm} 
\par
{\bf Proof:} Let $X \in \mfS_{g,b}^2(\CC)$. Stein factorisation
gives $X \to X' \to Y$, where $X'$ is a double covering of $Y$.
$X'$ has only rational double points
(and hence coincides with $X$), if and only if
$D$ has only simple singularities (\cite{BPV84} Th.\ II.5.1). 
Thus we have to show, that the open
subfunctor 
$$\widetilde{\mfS_{g,b}^2}(T) = \{ X \to T \in \mfS_{g,b}^2(T) | X \to 
\alpha(X) \hbox{ is finite} \}  $$
is also closed. Let $X \to T \in \mfS_{g,b}^2(T)$ be a family over the
pointed scheme $(T,0)$, such that
the restriction to $T' = T\sms 0$ is in $\widetilde{\mfS_{g,b}^2}(T')$.
After an {\'e}tale base change we may suppose, that $X$ admits a
simultaneous resolution of singularities. $X|_{T'}$ comes with
an involution $\tau$ whose quotient is $\alpha(X)|_{T'}$. As
$X$ has an ample canonical divisor, we can apply \cite{FP97}
Prop.\ 4.4 to extend $\tau$ to an involution on $X$. All we
need is to show, that for the fibre over $0$ we have $X_0/\tau \cong
\alpha(X_0)$. But 
this follows noting that $\alpha$ factors via $X/\tau$ and
that $X/\tau \to T$ is flat (see e.g.\ \cite{Ma97}).
\newline
The second assertion follows from the first and the theorem of 
Bertini, once we have shown that $L$ is ample. This immediately follows from
the exclusion of the case of surfaces isogenous to a product. 
And the last assertion should be clear by the given description.
\hfill $\Box$
\par
\begin{Lemma}
If $X \to T \in \mfS_{g,b}^2(T)$ admits a section, 
there are line bundles ${\cal L}_i$
on $B_i$, such that ${\cal L} = p_1^* {\cal L}_1 \otimes
p_2^* {\cal L}_2$.
\end{Lemma}
\par
{\bf Proof:} Due to the section we have
$$\PicF_{B_1\times B_2}/T \cong \PicF_{B_1/T} \times \PicF_{B_2} 
\times \Corr(B_1,B_2), $$
where $\Corr(B_1, B_2)$ is the functor of divisorial correspondences between
$B_1$ and $B_2$. If $T = \Spc \CC$ and $X$ is isotrivially fibred, we
know that ${L}^2 = {\cal O}_Y(D) = p_1^* 
{\cal O}_{B_1}(D_1) \otimes p_2^* {\cal O}_{B_2}(D_2)$. 
As $\Corr(B_1, B_2)$ has no nontrivial nilpotent
elements, this splitting is also possible for ${L}$ 
and as the zero section of $\Corr(B_1, B_2)$ is an open and
closed immersion, this remains also valid for arbitrary
$X \to T \in \mfS_{g,b}^2(T)$.
\hfill $\Box$
\par
By an {\'e}tale base change we can always suppose, that a section
exists. To analyse dimensions and connected components of the
moduli space, we can
hence use in the sequel the functor
$$Q(T) = \{ (Y/T,{\cal L}_1, 
{\cal L}_2, s)\} $$
where
$$ Y = B_1 \times B_2 \in \mfC_{b_1}(T) \times \mfC_{b_2}(T), 
{\cal L}_i \in \PicF_{{B}_i/T}, s \in \Gamma(Y, 
p_1^* {\cal L}_1 \otimes p_2^* {\cal L}_2) $$
instead of $\mfS_{g_1,b_1}^2(T)$.
\par
If we let $d_i = \deg {L}_i$, our description yields:
\par
\begin{Cor} \label{corcomp}
 Two surfaces in $\mfS_{g_1,b_1}^2(\CC)$ lie in the
same connected components, if and only if their invariants
$b_1$, $b_2$, $d_1$ and $d_2$ coincide. We hence denote these
components by $N_2(b_1,b_2,d_1, d_2)$.
\end{Cor}
\par
{\bf Proof:}
By the above theorem the conditions are clearly necessary.
\newline
Let ${\cal Y} = {\cal B}_1 \times {\cal B}_2 \to M^{[n]}$ denote
the universal family of products of smooth curves together with
level-$n$-structures (to obtain a fine moduli space).
Fix a line bundle $L = h_1^* L_1 \otimes h_2^* L_2$ on $Y$ and
let $g: P = \PicS_{{\cal Y}/M^{[n]}}^L \to M^{[n]}$ the
component of the relative Picard scheme parametrizing line
bundles linearly equivalent to $L$. Next we take the
scheme $F$ representing the functor
$$  T \mapsto \{({\cal Y}, {\cal L}, s) \,|\, ({\cal Y}, {\cal L}) \in
h_P(T), s \in \Gamma({\cal Y}, {\cal L}) \},$$
where $h_P$ denotes the functor of points of $P$. The existence of $F$
is garanteed by ([EGA] III, Th.\ 7.7.6) and
([EGA] I, Prop.\ 9.4.9). We will restrict $F$ to the open
locus $F'$, which parametrizes sections $s$ whose zero
locus has only simple singularities and we denote by $f$ the natural
morphism $F \to P$. The fibres of $f$ are
$H^0(C_1, L_1) \otimes H^0(C_2, L_2)$, where $C_i$ and $L_i$
are the curves and bundles corresponding to the image. And
of course, if this vector space is non-zero, the intersection with $F'$
is dense in it. 
\newline
We claim that $f(F')$ is connected: This is the
locus of quadruples $(C_i, L_i)$, $(i=1,2)$, such that 
$H^0(C_1, L_1) \otimes H^0(C_2, L_2)$ is non-zero.
$g$ maps $f(F')$ properly onto $M^{[n]}$, which is connected.
And the fibre of $g$ over $C_1 \times C_2$ is
$W^0_{d_1}(C_1) \times W^0_{d_2}(C_2)$, using the notation
in \cite{ACGH85}. This space is connected by \cite{ACGH85}
Th.\ V.1.4 and our claim follows by elementary topology.
\newline
Suppose $F'$ is the disjoint union of closed subsets
$A$ and $B$. Then $\ol{f(A)}$ and $\ol{f(B)}$ have a common point
$p \in f(F')$.  Suppose $p \in f(A)$ 
and hence $f^{-1}(p) \subset A$. Let $B_0 \subset
f(B)$ a subset, over which the fibres of $f$ have constant
positive dimension (i.e.\ $f|_{B_0}$ is a bundle of vector spaces) 
and such that $p \in \ol{B_0}$. Consider the closure
$C= \ol{f^{-1}(B_0)}$ in $F$, which is still
a bundle of vector spaces and whose image under $f$ contains $p$.
$C \cap F'$ is dense in each fibre of $f|_C$ 
and by construction $B \cap C$ is dense in $C$. This 
implies that $f^{-1}(p) \cap B \neq \emptyset$, a contradiction.
\hfill $\Box$
\par
For large $d_i$ we can easily compute the dimension of the
moduli space and the locus. Let
$h^0_i = \dim H^0(B_i, L_i^{\otimes 2})$ for $i=1,2$.
\par
\begin{Cor} \label{cordim}
If $2d_i > 2b_i -2$ (i.e $L^{\otimes 2}-K_Y$ is ample), the
dimension of $N_2(b_1,b_2,d_1, d_2)$ is
$$ \dim N_2 = 4b_1 + 4 b_2 - 7 + h^0_1 h^0_2$$
\newline
Isotrivial surfaces form irreducible, closed subvarieties in 
$N_2(b_1,b_2,d_1, d_2)$ of dimension
$ 4b_1 + 4b_2 + h^0_1$ and of dimension $4b_1 +4b_2 +h^0_2$.
\end{Cor}
\par
{\bf Proof:} $N_2(b_1,b_2,d_1, d_2)$ has a canonical
morphism to the moduli space of pairs of curves together with
a line bundle of degree $d_1$ and $d_2$ respectively. This
space is irreducible and has dimension $4b_1 + 4b_2-6$. The
ampleness of  $L^{\otimes 2}-K_Y$ ensures that
the fibres are dense in a vector space of dimension
$h^0_1 h^0_2 -1$ and this gives the first assertion.
\newline
We obtain the irreducible components of the locus of isotrivial
surfaces by taking only sections of the form $s=s^1 \otimes s^2$, 
i.e.\ by fixing $s^1$ and letting $s^2$ vary or vice versa.
The degree hypothesis implies that there are no obstructions
to deforming these sections.
\hfill $\Box$
\par
\begin{Rem} {\rm These results illustrates in the 'simplest possible' case 
Catanese's result (\cite{Ca00} Th.\ 5.4), 
that if we fix the 'topological type' (see loc.\ cit.\ for
details) of
the isotrivially fibred surface minus the singular fibres, 
we obtain irreducible components of the isotrivial locus, 
}
\end{Rem}
\par
%\begin{Rem}{\rm 
%
%}\end{Rem}
%\par

\section{The monodromy group is not invariant under deformations}

This section contains the 'negative' result, that in general one cannot
hope to analyse components of the moduli space containing isotrivial
surfaces by fixing the monodromy group.
\par
\begin{Prop} \label{propneq}
If the isotrivial surfaces $X_1$ and $X_2$ belong to one
connected component of $N_{g,b}$, the orders of the
respective monodromy groups $G_1$ and $G_2$ are equal, but
these groups need not to be isomorphic.
\end{Prop}
\par
{\bf Proof:} The first statement is clear because of Lemma
\ref{AlbLemma} and the fact that the group order is the
(generic) degree of the Albanese map.
\newline
To prove the second we can take the simplest possible case
$G_1 = \ZZ/4\ZZ$ and $G_2 = (\ZZ/2\ZZ)^2$ and use the
description of the isotrivial surfaces as abelian coverings
of $Y = B_1 \times B_2$. Fix $Y$ with $b_i \geq 2$ and
take line bundles $L_i$ on $B_i$ together with sections 
$s^i \in \Gamma(B_i, L_i^{\otimes 4})$. Let
$L = p_1^* L_1 \otimes p_2^* L_2$ and $s = s^1 \otimes s^2$.
The surface
$$ X_{G_1} = V(w_1^4 - s) \subset V(L^{\otimes 4}),$$
where $w_1$ is a local coordinate of $L$ and $V$ denotes the 
vanishing locus, is an isotrivial $G_1$-covering of $Y$. 
As $X_{G_1} = (C_1 \times C_2)/G_1$, where $C_i$ is
the $G_1$-covering corresponding to $L_i$ and $s^i$, this surface
has at most $A_3$ singularities. Due to the conditions
on $b_i$, no rational curve is contained in $X_{G_1}$, 
which is hence in $\mfS_{g,b}$ (where, as usual, $b=b_1$ and
$g$ is the genus of the fibres of $X \to B_1$).
\newline
Fix two generators $\alpha$ and $\beta$ of $G_2$ and two sections
$s_\alpha= s_\alpha^1 \otimes s_\alpha^2 \in \Gamma(Y, L^{\otimes 4})$,
$s_\beta= s_\beta^1 \otimes s_\beta^2 \in \Gamma(Y, L^{\otimes 2})$.
We can arrange that the zeros of $s_\alpha^i$ and $s_\beta^i$ ($i=1,2$)
are disjoint. Hence the surface
$$ X_{G_2} = V(w_\alpha^2 - s_\alpha, w_\beta^2 - s_\beta) \subset \VV(L^{\otimes 4}
\oplus L^{\otimes 2}),$$
where $w_\alpha$ and $w_\beta$ are local coordinates of $L^{\otimes 2}$
and $L$, has only $A_1$-singularities, i.e.\ is indeed in
$\mfS_{g,b}$.
\par
We construct a surface $X_1$ which is easily seen to deform to
both $X_{G_1}$ and $X_{G_2}$. Suppose there are
$s_2 \in \Gamma(Y,L^{\otimes 2})$ and $s_4 \in \Gamma(Y,L)$, 
such that $\div(s_2)$ and $\div(s_4)$ are both smooth and
have only normal crossings and such that
$\div(s_4+s_2^2/4)$ is smooth. We claim that, denoting
$w_i$ local coordinates of $L^{\otimes i}$,
$$X_1 = V(w_2^2 - s_2 w_2 -s_4, w_1^2-w_2) \subset \VV(L^{\otimes 4} 
\oplus L^{\otimes 2})$$
is in $\mfS_{g,b}$: The quotient $\tilde{X_1}$ of $X_1$ under 
$w_1 \mapsto -w_1$
ramifies over $\div(s_4+s_2^2/4)$ and is hence smooth. And
$X_1 \to \tilde{X_1}$ ramifies over $w_2=0$, which has smooth components
and normal crossings by hypotheses.
\newline
Obviously $X_1$ deforms to $V(w_1^4 -s_4) \subset \VV(L^{\otimes 4})$ and
this is a deformation of $X_{G_1}$. And by change of coordinates
$w_\alpha = w_2-s_2/2$, $w_1 = w_\beta$, we have
$$X_1 \cong  V(w_\alpha^2 - (s_4 + s_2^2/4), w_\beta^2 - w_\alpha - s_2/2)  
\subset  \VV(L^{\otimes 4} \oplus L^{\otimes 2}), $$
which obviously deforms to $X_{G_2}$.
\hfill $\Box$
\par

Martin M{\"o}ller \newline
Universit{\"a}t Karlsruhe \newline 
Mathematisches Institut II \newline
Englerstr. 2 \newline
D-76128 Karlsruhe, Germany \newline
e-mail: moeller@math.uni-karlsruhe.de \newline


\begin{thebibliography}{ABCD99}
\bibitem[ACGH85]{ACGH85} Arbarello, E., Cornalba, M., Griffiths, P.A., 
Harris, J., {\em Geometry of Algebraic Curves}, Volume I, Springer Grundlehren
267 (1985)
\bibitem[Be91]{Be91} Beauville, A., Letter to F.\ Catanese, Appendix
to [Ca91]
\bibitem[BPV84]{BPV84} Barth, W., Peters, C., van de Ven, A., {\em
Compact Complex Surfaces}, Ergebnisse der Math.\ 3 Bd.\ 4, 
Springer-Verlag (1984)
\bibitem[Ca91]{Ca91} Catanese, F., 
{\em Moduli and classification of irregular Kaehler
manifolds (and algebraic varieties) with Albanese general type fibrations}, 
Invent.\ Math.\ 104 (1991), 263--289
\bibitem[Ca00]{Ca00} Catanese, F., {\em Fibred surfaces, varieties
isogenous to a product and related moduli spaces}, Am. J. Math. 122 (2000),
1--44
\bibitem[FP97]{FP97} Fantechi, B., Pardini, R., {\em Automorphisms and
moduli spaces of varieties with ample canonical class via deformations
of abelian covers}, Comm.\ in Alg.\ 25(5) (1997), 1413--1441
\bibitem[Ma97]{Ma97} Manetti, M., {\em Iterated double covers and 
connected components of moduli spaces},
Topology 36 No.3 (1997), 745--764
\bibitem[Mo02]{Mo02}M{\"o}ller, M., {\em Modulr{\"a}ume irregul{\"a}r gefaserter
Fl{\"a}chen}, Dissertation, Karlsruhe (2002)
\bibitem[Ser96]{Ser96} Serrano, F., {\em Isotrivial Fibred Surfaces}, 
Annali di Matematica pura ed app.\ (IV), Vol.\ CLXXI (1996), 63--82
\end{thebibliography}
\end{document}